\documentclass[12pt,a4paper]{amsart}

\usepackage{amsfonts, amsmath, amssymb, amsthm, amscd}
\usepackage{anysize}
\usepackage{graphicx}

\usepackage{hyperref}
\hypersetup{
  colorlinks   = true, 
  urlcolor     = blue, 
  linkcolor    = blue, 
  citecolor   = red 
}

\usepackage{color}

\marginsize{2cm}{2cm}{2cm}{2cm}

\newtheorem{theorem}{Theorem}

\newtheorem*{claim*}{Claim}

\newtheorem{corollary}[theorem]{Corollary}

\theoremstyle{definition}

\theoremstyle{remark}
\newtheorem{remark}{Remark}

\newcommand*{\vv}[1]{ \overrightarrow{#1}}

\newcounter{fig}
\newcommand{\f}{\refstepcounter{fig} Fig. \arabic{fig}. }

\title{When different norms lead to same billiard trajectories?}
\author{Arseniy Akopyan}
\address{Arseniy Akopyan, Institute of Science and Technology Austria (IST Austria), Am Campus~1, 3400 Klosterneuburg, Austria}
\email{akopjan@gmail.com}

\author{Roman~Karasev}

\address{Roman Karasev,  Moscow Institute of Physics and Technology, Institutskiy per. 9, Dolgoprudny, Russia 141700 and Institute for Information Transmission Problems RAS, Bolshoy Karetny per. 19, Moscow, Russia 127994}
\email{r\_n\_karasev@mail.ru}
\urladdr{http://www.rkarasev.ru/en/}

\begin{document}


\begin{abstract}
In this paper, extending the works of Milena Radnovi\'{c} and Serge Tabachnikov, we establish conditions for two different non-symmetric norms to define the same billiard reflection law.
\end{abstract}
\maketitle

Milena Radnovi\'{c} in \cite{radnovic2003note} and independently Serge Tabchnikov in~\cite[Section 2]{tabachnikov2004remarks} made the following remarkable observation: 

\begin{theorem}
	\label{thm:radnovic-billiard}
	Let $\|\cdot\|_{\xi}$ be a not necessarily symmetric norm in the plane, having an ellipse with focus at the origin $o$ as the unit circle.
	Then $\|\cdot\|_{\xi}$ defines the same law of reflection as in the Euclidean metric: The angle of reflection equals to the angle of incidence.
\end{theorem}

In the first paper it was also notices that it leads to the fact that billiard trajectories in the plane with norm, defined by an ellipse as the unit circle, are the same as in the Euclidean plane after a suitably chosen affine transform.

Another consequence of their theorem is that the Euclidean and normed ellipses with foci $o$ and $f$ coincide, where $f$ is the second focus of $\xi$.
Indeed, if $x$ is a point in the plane, then by Theorem \ref{thm:radnovic-billiard} the differential of $\|\vv{ox}\|_{\xi} + \|\vv{xf}\|_{\xi}$ is proportional to the differential of the same expression for the Euclidean norm for every $x$. Hence the value $\|\vv{ox}\|_{\xi} + \|\vv{xf}\|_{\xi}$ does not change when $x$ moves along the ellipse with foci $f$ and $o$, along the zero direction of both differentials.

It is interesting, that the latter statement may be deciphered to the following elementary geometric formulation, for which we do not know any short synthetic proof essentially different from the one stated in the above paragraph:

\begin{corollary}
	\label{cor:eucliean ratios}
	Let $\xi_1$ and $\xi_2$ bet two confocal ellipses with foci at $f_1$ and $f_2$.
	For each point $x$ on $\xi_1$, denote by $y_1$ and $y_2$ the points of intersections of $\xi_2$ with rays $f_1x$ and $f_2x$ respectively. Then for any point $x$ on $\xi_1$:
	\[
	\frac{|xf_1|}{|y_1f_1|}+\frac{|xf_2|}{|y_2f_2|} =\mathrm{const}.
	\]
\end{corollary}
It can be shown, that the constant in the corollary above equals:
$\frac{\ell_1-|f_1f_2|}{\ell_2-|f_1f_2|}+\frac{\ell_1+|f_1f_2|}{\ell_2+|f_1f_2|}$, where  $\ell_1$ and $\ell_2$ are the major axes of the ellipses $\xi_1$ and $\xi_2$.

\begin{center}
	\includegraphics{milena-fig-1.mps}\\
	\f \label{fig:two-ellipses} 
\end{center}

Now we extend the Radnovi\'{c}--Tabachnikov theorem to normed spaces in higher dimension:
\begin{theorem}
	\label{thm:main theorem}
	Let $K$ be a smooth convex body in $\mathbb{R}^n$ containing the origin, and $T$ be its convex image under a projective transform, which maps each line passing through origin to itself preserving its orientation at the origin.
	Then the billiard reflection law in the space with norm $\|\cdot\|_{K}$ is the same as in the space with norm $\|\cdot\|_{T}$.
\end{theorem}

\begin{remark}
	It is known (see \cite[Lemma 4.6]{akopyan2018incircular}) that such kind of projective transforms send spheres with center at origin to ellipsoids of rotation with one of the foci at the origin.
	Therefore this theorem directly generalizes Theorem~\ref{thm:radnovic-billiard} to higher dimension.
\end{remark}

\begin{remark}
	The law of reflection is not well-defined for convex bodies $K$, which are not strictly convex.
	In this case we may follow the conventions in \cite{akopyan2015billiards} and define billiard trajectories, for which the reflection direction is not uniquely defined. 
\end{remark}

\begin{center}
	\includegraphics{milena-fig-2.mps}\\
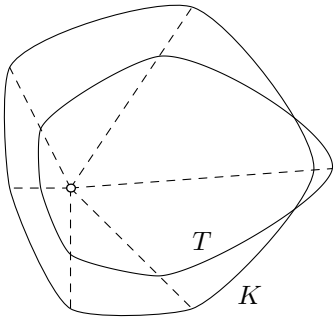
	\f \label{fig:body and image } Convex body $K$ and its projective image $T$.
\end{center}

\begin{proof}
	We use ideas from \cite{gutkin2002billiards, artstein2008brunn, akopyan2014elementary}, suggesting to work with billiard trajectory in the Banach space $U=\mathbb{R}^n$ with norm $\|\cdot\|_{K}$ in terms of momenta in the dual space $U^*$ with norm $\|\cdot\|_{K^\circ}$, whose unit ball is the polar body~$K^\circ$.
	From the smoothness of $K$, to each unit velocity $u\in \partial K$ there corresponds a conjugate unit momentum $u^* \in \partial K^\circ$, such that $u^*(u) = 1$ and $\|u^*\|_{K^\circ} = 1$. The equation $u^*(x) = 1$ defines the hyperplane tangent to $\partial K$ at $u$, while the equation $u(y) = 1$  defines the hyperplane tangent to $\partial K^\circ$ at $u^*$.
	
	Let $q_1q_2q_3$ be a part of a billiard trajectory in $U$, where $q_2$ is the point where the trajectory hits a hypersurface $S$ and reflects.
	Then the sum $\|\vv{q_1x}\|_K+\|\vv{xq_3}\|_K$ as a function of $x \in S$ has a critical value at $q_2$.
	The criticality in terms of first derivatives means:
\begin{equation}
		u_2^*-u_1^*=\lambda n^*,
		\tag{$*$} \label{eq:dual reflection law}
\end{equation}
	where $u_1^*$ and $u_2^*$ are momenta corresponding to unit vectors in directions $\vv{q_1q_2}$ and $\vv{q_2q_3}$, and $n^*$ is the normal covector to $S$ at $q_2$.
	
	It is crucial that equation \eqref{eq:dual reflection law} is preserved under a positive similiarity of the body $K^\circ$, possibly with different factor $\lambda$.
	Indeed, let $T^\circ=tK^\circ+v^*$, where $t>0$ and $v^*\in U^*$.
	It is easy to see that the momentum, corresponding to velocity $u$ with respect to the body $T^\circ$, equals $u_T^*=tu^*+v^*$, because $u$ is a linear function on $U^*$ and the points, where its maximum is obtained on $K^\circ$ and $T^\circ$ are moved one to another by the homothety.
	Hence the difference of the new momenta at a reflection point equals to $t(u_2^*-u_1^*)$, which is still parallel to $n^*$.
	
	We obtain that the reflection laws for two norms $\|\cdot\|_{K}$ and $\|\cdot\|_{T}$ coincide if $K^\circ$ and $T^\circ$ are positive homothets of each other.  
	A positive homothety is a projective transform which maps any point at infinity to itself.
	Therefore in the dual space (our original $U$) a positive homothety corresponds to the map, which preserves its polar images as a sets, that is planes passing through the origin.
	It is easy to see, that this is the projective transform described in the statement of the theorem.
	
	In simpler words, $K$ is given by the system of linear inequalities of the form
	\[
	u^*(x) \le 1,\quad \forall u^*\in K^\circ.
	\]
	Hence the equations of $T$ must be (assuming working not far from the origin, where $v(x) < 1$)
	\[
	t u^*(x) + v^*(x) \le 1\Leftrightarrow u^*\left( \frac{tx}{1 - v(x)}\right) \le 1, \quad \forall u^*\in K^\circ.
	\]
	It remains to note that 
	\[
	x \mapsto \frac{tx}{1 - v(x)}, \quad t > 0
	\]
	is the general form of projective maps that preserve lines thorough the origin and keep their orientations at the origin.

		%
	%
\end{proof}

\section*{Acknowledgments}
The authors thank Alexey Balitskiy, Milena Radnovi\'{c}, and Serge Tabachnikov for useful discussions.

AA was supported by European Research Council (ERC) under the European Union's Horizon 2020 research and innovation programme (grant agreement No 78818 Alpha).
RK was supported by the Federal professorship program grant 1.456.2016/1.4 and the Russian Foundation for Basic Research grants 18-01-00036 and 19-01-00169.

\bibliographystyle{abbrv}
\bibliography{billiards.bib}

\begin{thebibliography}{1}

\bibitem{akopyan2015billiards}
A.~Akopyan and A.~Balitskiy.
\newblock Billiards in convex bodies with acute angles.
\newblock {\em Israel Journal of Mathematics}, 216(2):833--845, 2016.

\bibitem{akopyan2014elementary}
A.~Akopyan, A.~Balitskiy, R.~Karasev, and A.~Sharipova.
\newblock Elementary approach to closed billiard trajectories in asymmetric
  normed spaces.
\newblock {\em Proc. Amer. Math. Soc.}, 144(10):4501--4513, 2016.

\bibitem{akopyan2018incircular}
A.~V. Akopyan and A.~I. Bobenko.
\newblock Incircular nets and confocal conics.
\newblock {\em Transactions of the American Mathematical Society},
  370(4):2825--2854, nov 2018.

\bibitem{artstein2008brunn}
S.~Artstein-Avidan and Y.~Ostrover.
\newblock A {B}runn--{M}inkowski inequality for symplectic capacities of convex
  domains.
\newblock {\em IMRN: International Mathematics Research Notices}, 2008, 2008.

\bibitem{gutkin2002billiards}
E.~Gutkin and S.~Tabachnikov.
\newblock Billiards in {F}insler and {M}inkowski geometries.
\newblock {\em J. Geom. Phys.}, 40(3-4):277--301, 2002.

\bibitem{radnovic2003note}
M.~Radnovi\'{c}.
\newblock A note on billiard systems in {F}insler plane with elliptic
  indicatrices.
\newblock {\em Publ. Inst. Math. (Beograd) (N.S.)}, 74(88):97--101, 2003.

\bibitem{tabachnikov2004remarks}
S.~Tabachnikov.
\newblock Remarks on magnetic flows and magnetic billiards, {F}insler metrics
  and a magnetic analog of {H}ilbert's fourth problem.
\newblock In {\em Modern dynamical systems and applications}, pages 233--250.
  Cambridge Univ. Press, Cambridge, 2004.

\end{thebibliography}

\end{document}